\documentclass{gtmon_a}
\pdfoutput=1
\usepackage{amssymb} 
\usepackage{amsmath,amscd} 

\proceedingstitle{The interaction of finite-type and Gromov--Witten
invariants (BIRS 2003)}
\conferencestart{15 November 2003}
\conferenceend{20 November 2003}
\conferencename{The interaction of finite-type and Gromov--Witten
invariants}
\conferencelocation{Banff International Research Station, Banff, Alberta,
Canada}

\editor{David Auckly}
\givenname{David}
\surname{Auckly}

\editor{Jim Bryan}
\givenname{Jim}
\surname{Bryan}

\title{Perturbative expansion of Chern--Simons theory}                    

\author{Justin Sawon}
\givenname{Justin}
\surname{Sawon}
\address{Department of Mathematics\\
SUNY at Stony Brook\\\newline
NY 11794\\USA}
\email{sawon@math.sunysb.edu}
\urladdr{http://www.math.sunysb.edu/~sawon/}

\volumenumber{8}
\issuenumber{}
\publicationyear{2006}
\papernumber{7}
\startpage{145}
\endpage{166}

\doi{}
\MR{}
\Zbl{}


\received{10 December 2004}
\revised{25 April 2005}
\accepted{14 December 2004}
\published{22 April 2006}
\publishedonline{22 April 2006}
\proposed{}
\seconded{}
\corresponding{}
\editor{}
\version{}

\arxivreference{math.GT/0504495}

\begin{asciiabstract}
An overview of the perturbative expansion of the Chern--Simons 
path integral is given. The main goal is to describe how trivalent graphs
appear: as they already occur in the perturbative expansion of an
analogous finite-dimensional integral, we discuss this case in
detail.
\end{asciiabstract}

\makeatletter
\def\cnewtheorem#1[#2]#3{\newtheorem{#1}{#3}[section]
\expandafter\let\csname c@#1\endcsname\c@thm}

\AtBeginDocument{\let\bar\wbar\let\tilde\wtilde\let\hat\what}

\cnewtheorem{lem}[thm]{Lemma}
\cnewtheorem{conj}[thm]{Conjecture}

\theoremstyle{definition}
\cnewtheorem{defn}[thm]{Definition}
\newtheorem*{rem}{Remark}
\newtheorem*{exm}{Example}

\makeatother  

\newcommand{\thetagraph}{\raisebox{0pt}{
	\begin{picture}(20,15)(-10,-5)
	\put(0,0){\circle{20}} \put(-10,0){\line(1,0){20}}
	\end{picture}}}

\newcommand{\dumbell}{\raisebox{0pt}{
        \begin{picture}(42,15)(-21,-5)
        \put(-14,0){\circle{14}} \put(14,0){\circle{14}} 
  	\put(-7,0){\line(1,0){14}}
	\end{picture}}}

\newcommand{\bigdumbell}{\raisebox{0pt}{\small
	\begin{picture}(200,60)(-100,-30)
	\put(-60,0){\circle{40}} \put(60,0){\circle{40}}
	\put(-40,0){\line(1,0){80}}
	\put(-62,-3){$V_{jkl}$} \put(41,-3){$V_{abc}$}
	\put(-100,-3){$Q^{jk}$} \put(-5,3){$Q^{la}$}
	\put(81,-3){$Q^{bc}$}
	\put(-40,8){$\scriptstyle j$} \put(-40,-11){$\scriptstyle k$} 
        \put(-31,2){$\scriptstyle l$} \put(35,8){$\scriptstyle b$} 
        \put(34,-11){$\scriptstyle c$} \put(26,2){$\scriptstyle a$}
	\end{picture}}}

\newcommand{\involutionI}{\raisebox{0pt}{
        \begin{picture}(42,15)(-21,-5)
        \put(-14,0){\circle{14}} \put(14,0){\circle{14}} 
  	\put(-7,0){\line(1,0){14}} \put(-18,-2.5){$\updownarrow$}
	\end{picture}}}

\newcommand{\involutionII}{\raisebox{0pt}{
        \begin{picture}(42,15)(-21,-5)
        \put(-14,0){\circle{14}} \put(14,0){\circle{14}} 
  	\put(-7,0){\line(1,0){14}} \put(-6.5,1){$\leftrightarrow$}
	\end{picture}}}

\newcommand{\involutionIII}{\raisebox{0pt}{
        \begin{picture}(42,15)(-21,-5)
        \put(-14,0){\circle{14}} \put(14,0){\circle{14}} 
  	\put(-7,0){\line(1,0){14}} \put(10,-2.5){$\updownarrow$}
	\end{picture}}}

\newcommand{\involutionIV}{\raisebox{0pt}{
	\begin{picture}(20,15)(-10,-5)
	\put(0,0){\circle{20}} \put(-10,0){\line(1,0){20}}
	\put(-6.5,1){$\leftrightarrow$}
	\end{picture}}}

\newcommand{\vertex}{\raisebox{0pt}{
	\begin{picture}(20,15)(-10,-3)
	\put(-8,0){\line(1,0){16}} \put(-8,0){\line(2,1){16}}
	\put(-8,0){\line(2,-1){16}}
	\end{picture}}}

\newcommand{\edge}{\raisebox{0pt}{
	\begin{picture}(20,15)(-10,-3)
	\put(-8,0){\oval(18,12)[r]}
	\end{picture}}}

\newcommand{\Tr}{\operatorname{Tr}}
\newcommand{\Aut}{\operatorname{Aut}}
\newcommand{\Vol}{\operatorname{Vol}}
\newcommand{\overc}{{}\mskip1mu\overline{\mskip-1mu c \mskip-1mu}\mskip1mu}

\begin{document}

\begin{abstract} 
An overview of the perturbative expansion of the Chern--Simons 
path integral is given. The main goal is to describe how trivalent graphs
appear: as they already occur in the perturbative expansion of an
analogous finite-dimensional integral, we discuss this case in
detail.
\end{abstract}

\maketitle

\section{Introduction}

While a standard technique used by physicists, many
mathematicians are not familiar with the mechanisms and justifications
behind perturbatively expanding Feynman path integrals. In this
article we will describe how trivalent graphs arise when one expands a
path integral whose Lagrangian contains a cubic term. We'll start with
a finite-dimensional integral, and then indicate the necessary
adjustments required to generalize to the infinite-dimensional
integral arising in Chern--Simons theory.

This article is for the proceedings of the BIRS workshop ``The
interaction of finite type and Gromov--Witten invariants''. The purpose
of the workshop was to investigate `large--$N$ duality', which relates 
perturbative Chern--Simons theory with gauge group $\mathrm{SU}(N)$ to
certain open and closed string theories. One of the more mysterious
aspects of this correspondence, to this author at least, is that in
string theory trivalent graphs (more precisely, fat graphs) have some 
geometric realization, whereas in perturbative Chern--Simons theory
they enter purely as part of the technique of calculation. As we shall
see, trivalent graphs arise merely as combinatorial objects which
index terms in the series expansion.

The exact solution to Chern--Simons theory was described by
Witten~\cite{wittenEightNine} in the late 1980s. The perturbative expansion of
the path integral for flat $\mathbb{R}^3$ was then considered by
Guadagnini, Martellini, and 
Mintchev~\cite{gmmEightNineiii,gmmEightNineii,gmmEightNinei,gmmNineZero}, and for general
three-manifolds by Axelrod and
Singer~\cite{asNineTwo,asNineFour}. Subsequently much attention was focused on
understanding the perturbative series for knots and links in
$\mathbb{R}^3$ (see Bar-Natan~\cite{bar-natanNineOne,bar-natanNineFive}, Bott and
Taubes~\cite{btNineFour}, Altschuler and Freidel~\cite{afNineSeven}, and
Thurston~\cite{thurstonNineFive}). More on perturbative expansion of
path integrals may be found in Witten's lectures~\cite{wittenNineNine}. For
a slightly different application, to the perturbative expansion of
Hermitian matrix integrals, see Mulase~\cite{mulaseNineFour}, where Feynman
diagram techniques are nicely explained.

While this article was being prepared, a preprint by Michael
Polyak~\cite{polyakZeroFour} appeared, covering very similar material. All
the same, the author decided to complete this article for inclusion in
the present volume. The interested reader may wish to consult Polyak's
preprint, which in many cases provides additional details to what is
covered here.

\subsection*{Acknowledgements}

The author is very grateful to Alexander Kirillov Jr for reading an
earlier version of this article and suggesting numerous improvements. 
The referee also made many helpful comments. The author would like to
thank the organizers of the workshop for putting together such an
interesting programme, and MSRI for travel support. The author is
partially supported by NSF grant 0305865.

\section{Stationary-phase approximation}

Fix a compact, connected, simply-connected gauge group $G$ (for
example, $\mathrm{SU}(n)$). These conditions imply that
$\pi_1(G)=\pi_2(G)=1$ and $\pi_3(G)=\mathbb{Z}$, a fact that we shall
use shortly. Let $M$ be a three-manifold without boundary. The field
in Chern--Simons theory is a connection $\nabla$ on a principal
$G$--bundle over $M$. The conditions on $G$ imply that the bundle can
be trivialized and hence we can write the connection as $\nabla=d+A$
where $A\in\Omega^1(M)\otimes\mathfrak{g}$ is a Lie algebra valued
one-form. A change of trivialization is given by a gauge
transformation, which is a smooth function $g\co M\rightarrow G$. Under
this transformation $A$ changes to $g^{-1}Ag+g^{-1}dg$.

\begin{defn}
The Chern--Simons action is given by
$$S(A)=\frac{1}{4\pi}\int_M
  \Tr\bigl(A\wedge dA+{\textstyle\frac{2}{3}}A\wedge A\wedge A\bigr).$$
\end{defn}

\begin{rem}
Let us explain the terms in the integrand. Write $A=\zeta^IE_I$ where
$\zeta^I$ are one-forms and $\{E_I\}$ is a basis for $\mathfrak{g}$.
On forms, $\wedge$ denotes the usual wedge product. On the Lie algebra
part, the first term $\Tr(E_IE_J)$ is an invariant bilinear
form $B_{IJ}$ normalized so that $B(h_{\alpha},h_{\alpha})=2$, where
$h_{\alpha}$ is a coroot corresponding to a long root $\alpha$. The
second term looks like
\begin{eqnarray*}
\Tr(A\wedge A\wedge A) & = & \Tr(E_IE_JE_K)\zeta^I\!\!\wedge
\zeta^J\!\!\wedge \zeta^K \\
   & = & {\textstyle\frac{1}{2}}\Tr([E_I,E_J]E_K)\zeta^I\!\!\wedge
\zeta^J\!\!\wedge \zeta^K \\
   & = & {\textstyle\frac{1}{2}}c_{IJ}^LB_{LK}\zeta^I\!\!\wedge
\zeta^J\!\!\wedge \zeta^K
\end{eqnarray*}
where $c_{IJ}^K$ are the structure constants of $\mathfrak{g}$. We
usually write $c_{IJ}^LB_{LK}$ simply as $c_{IJK}$. For a matrix group
such as $\mathrm{SU}(n)$, these operations can be performed simply as
they appear: by multiplying and taking traces of matrices.

Observe that $c_{IJK}$ gives a three-form $\omega$ on $\mathfrak{g}$
defined by
$$\omega(X,Y,Z):=B([X,Y],Z).$$
where $X$, $Y$, $Z\in\mathfrak{g}$ and $B$ is the invariant bilinear
form as above. This can be made into a left-invariant three-form on
$G$, and $\int_{\gamma}\omega=8\pi^2$ where $\gamma$ is a generator of
$\mathrm{H}_3(G,\mathbb{Z})=\mathbb{Z}$.
\end{rem}

The action $S$ gives a functional on the space of connections
$\mathcal{A}$. Let us calculate its derivative.

\begin{lem}
\label{stationary=flat}
The curvature $F$ of $A$ can be regarded as a one-form on
$\mathcal{A}$, and the derivative of $S$ equals $\frac{1}{2\pi}F$.
\end{lem}

\begin{proof}
Let $\delta A$ describe a variation of the connection $A$. We can
regard $\delta A$ as an element of $T_A\mathcal{A}$, and then the
curvature $F$ of $A$ defines a one-form by
$$\delta A \mapsto \int_M\Tr(\delta A\wedge F).$$
Moreover,
\begin{eqnarray*}
S(A+\delta A)-S(A) & = & \frac{1}{4\pi}\int_M\Tr(\delta A\wedge
dA+A\wedge d\delta A+2\delta A\wedge A\wedge A) \\
 & = & \frac{1}{2\pi}\int_M\Tr(\delta A\wedge(dA+A\wedge A))
\\
 & = & \frac{1}{2\pi}\int_M\Tr(\delta A\wedge F)
\end{eqnarray*}
where we have used integration by parts.
\end{proof}

\begin{rem}
The gauge group $\mathcal{G}$ acts on $\mathcal{A}$, and $F$ is an
invariant one-form. On the other hand, a calculation shows that
$$S(g^{-1}Ag+g^{-1}dg)-S(A)=-2\pi\deg (g)$$
where $\deg (g)$ is the degree of $g$ as a map from $M$ to $G$, that is,
$g$ takes the fundamental class of $M$ to $\deg (g)\gamma$, where
$\gamma$ is the generator of $\mathrm{H}_3(G,\mathbb{Z})$ mentioned
earlier. Note that the degree gives an isomorphism
$\pi_0(\mathcal{G})\cong\mathbb{Z}$.
\end{rem}

Since $S$ is well-defined on $\mathcal{A}/\mathcal{G}$ up to a
multiple of $2\pi$, the integrand in the next definition is invariant
under gauge transformations.

\begin{defn}
The partition function of Chern--Simons theory is given by the path
integral
$$Z_k(M)=\int_{\mathcal{A}}\mathcal{D}Ae^{ikS(A)}$$
where the level $k$ is a positive integer.
\end{defn}

\begin{rem}
The measure on the infinite-dimensional space $\mathcal{A}$ of
connections is not well-defined, so the path integral is
really a heuristic tool. It suggests that there should be a
topological invariant $Z_k(M)$ of $M$. Moreover, the path integral
suggests various ways of studying this invariant, for instance via a
topological quantum field theory or via a perturbative series.
\end{rem}

In the perturbative approach $k$ should be thought of as being
proportional to the inverse of Planck's constant $\hbar^{-1}$. In the
classical limit $\hbar\rightarrow 0$, so we should expand in
$k^{-1}$. In the $k\rightarrow\infty$ limit the path integral
localizes around stationary points: this should be understood as
follows. Assume that $A$ is not a stationary point of the action
$S$. The integrand $\exp(ikS(A))$ lies on the unit circle, and if we
vary $A$ slightly the value of the integrand rotates around the
circle; in the large--$k$ limit the values that the integrand takes are
distributed evenly around the unit circle and so cancel out in the
integral.

By \fullref{stationary=flat} the derivative of $S$ vanishes at flat
connections, so these are the stationary points. To expand around a
flat connection $\alpha$, we let $A=\alpha+\beta$ where
$\beta\in\Omega^1(M)\otimes\mathfrak{g}$ should be regarded as being
the perturbation from $\alpha$. A calculation shows that
$$S(A)=S(\alpha)+\frac{1}{4\pi}\int_M\Tr\bigl(\beta\wedge
  d_{\alpha}\beta+{\textstyle\frac{2}{3}}\beta\wedge\beta\wedge\beta\bigr)$$
where $d_{\alpha}$ is the connection on one-forms induced by the
connection $\alpha$. The contribution to $Z_k(M)$ from the flat
connection $\alpha$ is therefore
\begin{eqnarray}
\label{infinite}
\int_{\beta\in\Omega^1(M)\otimes\mathfrak{g}}d\beta
e^{ikS(\alpha)}e^{\frac{ik}{4\pi}\int_M\Tr(\beta\wedge
d_{\alpha}\beta+\frac{2}{3}\beta\wedge\beta\wedge\beta)}
\end{eqnarray}
where the first factor $\exp(ikS(\alpha))$ in the integrand is a
constant and the second factor involves both a quadratic and a cubic
term in $\beta$. Note that the linear term in $\beta$ vanishes because
we are expanding around a stationary point.

This process of localization reduces the path integral to an integral
over the moduli space of flat connections. If the flat connections are
isolated then $Z_k(M)$ is simply a sum of integrals like the one
above. However, the flat connections are never isolated because of the
action of the gauge group - we will return to this point later. Note
also that \eqref{infinite} is still an infinite-dimensional integral;
our first step will be to calculate the integral in the simpler
finite-dimensional case.

\section{A finite-dimensional model}

We will consider the following finite-dimensional integral
$$Z_k=\int_{\mathbb{R}^N}d^Nxe^{ik\bigl(\frac{1}{2}Q(x,x)+\frac{1}{6}V(x,x,x)\bigr)}$$
where $Q$ and $V$ are respectively quadratic and cubic forms on
$\mathbb{R}^N$ (linear in all entries and totally symmetric), and $k$
is a positive real constant. One could drop the symmetry assumption on
$V$, which would require some additional labeling of the half-edges at
trivalent vertices of the graphs arising later, but in this article we
will only consider symmetric $V$, which is why we have included the
factor of $1/3!$. 

The integral $Z_k$ is not absolutely convergent, but it can be defined
in some formal sense via analytic continuation. To begin with, the
Gaussian integral on $\mathbb{R}$ can be evaluated
$$\int_{\mathbb{R}}dxe^{-\frac{1}{2}kx^2}=\sqrt{\frac{2\pi}{k}}$$
for $k$ real and positive. Analytic continuation in $k$ yields
$$\int_{\mathbb{R}}dxe^{\frac{1}{2}ikx^2}=\left\{\begin{array}{cc}
   \sqrt{\frac{2\pi}{k}}e^{\frac{\pi i}{4}} & \text{if}\ k>0, \\
   \sqrt{\frac{2\pi}{|k|}}e^{-\frac{\pi i}{4}} & \text{if}\ k<0.
   \end{array}\right.$$
In dimension greater than one, we can diagonalize the quadratic form
to get
$$\int_{\mathbb{R}^N}d^Nxe^{-k\frac{1}{2}Q(x,x)}=\Biggl(\frac{2\pi}{k}\Biggr)^{\frac{N}{2}}\frac{1}{(\det Q)^{\frac{1}{2}}}$$
for $Q$ positive definite, and
$$\int_{\mathbb{R}^N}d^Nxe^{ik\frac{1}{2}Q(x,x)}=\Biggl(\frac{2\pi}{k}\Biggr)^{\frac{N}{2}}\frac{e^{\pi
i(\mathrm{sign}Q)/4}}{\bigl|\det Q\bigr|^{\frac{1}{2}}}$$
for $Q$ indefinite (but still non-degenerate). We will refer to all of
these integrals, which don't involve cubic terms, as \emph{Gaussians}.

Returning to $Z_k$, we first deal with the dependence on $k$ by a
change of variable $x\mapsto x^{\prime}:=\sqrt{k}x$ (we immediately
rewrite the dummy variable $x^{\prime}$ as $x$ in our formula)
\begin{eqnarray*}
Z_k & = &
k^{-\frac{N}{2}}\int_{\mathbb{R}^N}d^Nxe^{i\bigr(\frac{1}{2}Q(x,x)+\frac{1}{6}k^{-\frac{1}{2}}V(x,x,x)\bigl)} \\
 & = &
k^{-\frac{N}{2}}\int_{\mathbb{R}^N}d^Nxe^{i\frac{1}{2}Q(x,x)}\sum_{m=0}^{\infty}\frac{i^m}{m!6^mk^{\frac{m}{2}}}V(x,x,x)^m.
\end{eqnarray*}
We have used the power series for the exponential function involving
the cubic term, but notice that the odd terms in the series give an
integrand which is an odd function on $\mathbb{R}^N$, so those terms
do not contribute to the integral. On the other hand, the
$2m$th term gives the coefficient of $k^{-\frac{N}{2}-m}$
as
$$\int_{\mathbb{R}^N}d^Nxe^{i\frac{1}{2}Q(x,x)}
  \frac{(-1)^m}{(2m)!6^{2m}}V(x,x,x)^{2m}.$$
To calculate this we need to understand how to integrate the Gaussian
multiplied by a polynomial.

\begin{lem}
\label{lemThreeOne}
Let $J(x)=J_jx^j$ be a linear functional on $\mathbb{R}^N$ (we are
using the Einstein summation convention) and define
$$Z(J):=\int_{\mathbb{R}^N}d^Nxe^{i\bigl(\frac{1}{2}Q(x,x)+J(x)\bigr)}$$
so that $Z(0)$ is just the Gaussian integral. Then
$$Z(J)=Z(0)e^{-\frac{i}{2}Q^{-1}(J,J)}$$
where $Q^{-1}(J,J)=Q^{jk}J_jJ_k$ and $Q^{jk}$ is the matrix of the
inverse of $Q$.
\end{lem}

\begin{proof}
This is a simple exercise in completing the square.
\end{proof}

The linear term $J$ is known as a \emph{source}, and the main trick is to
differentiate $Z(J)$ with respect to $J$ and then evaluate at
$J=0$. We will illustrate this through some examples.

\begin{exm}
First observe that, using the definition of $Z(J)$,
$$\int_{\mathbb{R}^N}d^Nxe^{i\frac{1}{2}Q(x,x)}x^j=
  \biggl(-i\frac{\partial\phantom{J_j}}{\partial
J_j}\biggr)Z(J)\bigr|_{J=0}.$$
Now using the expression for $Z(J)$ in \fullref{lemThreeOne}, the
right-hand side becomes
$$-iZ(0)\bigl(-iQ^{jk}J_ke^{-\frac{i}{2}Q^{-1}(J,J)}\bigr)\bigr|_{J=0}.$$
This vanishes when we evaluate at $J=0$, as it should because it is
the integral over $\mathbb{R}^N$ of an odd function.

Next observe that to integrate the product of the Gaussian and a
degree two polynomial, we should differentiate $Z(J)$ twice
$$\int_{\mathbb{R}^N}d^Nxe^{i\frac{1}{2}Q(x,x)}x^jx^k=\biggl(-i\frac{\partial\phantom{J_j}}{\partial
J_j}\biggr)\biggl(-i\frac{\partial\phantom{J_k}}{\partial
J_k}\biggr)Z(J)\bigr|_{J=0}.$$
Once again, we use \fullref{lemThreeOne} to evaluate this as
$$(-i)^2Z(0)
  \bigl(-iQ^{jk}e^{-\frac{i}{2}Q^{-1}(J,J)}+J_k(\cdots)\bigr)\bigr|_{J=0}
  = iZ(0)Q^{jk}.$$  
The precise expression $(\cdots)$ is irrelevant as it
vanishes when we evaluate at $J=0$.

In general, to integrate the product of the Gaussian and a degree $d$
polynomial, we should differentiate $Z(J)$ $d$ times. 
\end{exm}

Now let us apply this to the $m=1$ term occurring in $Z_k$
$$\int_{\mathbb{R}^N}d^Nxe^{i\frac{1}{2}Q(x,x)}\left(\frac{-1}{2!6^2}\right)V(x,x,x)^2.$$
Writing $V(x,x,x)$ in coordinates as $V_{jkl}x^jx^kx^l$ or
$V_{abc}x^ax^bx^c$, we see that the integrand is the product of a
Gaussian and a degree six polynomial, and hence the integral equals
\begin{multline*}
 \left(\frac{-1}{2!6^2}\right)V_{jkl}V_{abc}
  \left(-i\frac{\partial\phantom{J_j}}{\partial J_j}\right)
  \left(-i\frac{\partial\phantom{J_k}}{\partial J_k}\right)
  \cdots\left(-i\frac{\partial\phantom{J_c}}{\partial J_c}\right)
  Z(J)\biggr|_{J=0} \\
= \left(\frac{-1}{2!6^2}\right)
  V_{jkl}V_{abc}i^3Z(0)
  \bigl(Q^{jk}Q^{la}Q^{bc}+Q^{jk}Q^{lb}Q^{ac}+\cdots\bigr)
\end{multline*}
The last factor on the right consists of the sum over all ways of
dividing the set of indices $\{j,k,l,a,b,c\}$ into pairs, and hence
contains $\smash{\frac{1}{3!}\binom{6}{\,2,2,2\,}}=15$ terms. Any other terms
which appear when we differentiate $Z(J)$ will vanish when we evaluate
at $J=0$.

Since we are using the summation convention, $j,k$, etc, are dummy
indices. Up to a factor of $\bigl(\frac{i}{2!6^2}\bigr)Z(0)$, we get
nine terms of the form $V_{jkl}V_{abc}Q^{jk}Q^{la}Q^{bc}$ and six of
the form
$V_{jkl}V_{abc}Q^{ja}Q^{kb}Q^{lc}$. These terms can be encoded using
the trivalent graphs
$$\Gamma_1=\dumbell\qquad\mbox{and}\qquad\Gamma_2=\thetagraph.$$
More precisely, to each trivalent graph $\Gamma$ we associate a weight
$W(\Gamma)$ given by labeling vertices with $V$ and edges with
$Q^{-1}$, and then contracting indices. For instance on $\Gamma_1$ we
first label the legs (ends of the edges) with indices. Then we label
the vertices with $V_{jkl}$ and $V_{abc}$, and the edges with
$Q^{jk}$, $Q^{la}$, and $Q^{bc}$, in such a way that the indices agree
with those on the legs.
$$\bigdumbell$$
Then we apply the summation convention and arrive at the weight
\begin{eqnarray}
\label{summation_convention}
W(\Gamma_1) & := & V_{jkl}V_{abc}Q^{jk}Q^{la}Q^{bc}.
\end{eqnarray}
In this construction $V$ and $Q^{-1}$ are known as the \emph{vertex} and
\emph{propagator}, respectively.

Summing up, the $m=1$ term becomes
\begin{eqnarray*}
\biggl(\frac{i}{2!6^2}\biggr)Z(0)\bigl(9W(\Gamma_1)+6W(\Gamma_2)\bigr)
 & = & iZ(0)\bigl({\textstyle\frac{1}{8}}W(\Gamma_1)
   +   {\textstyle\frac{1}{12}}W(\Gamma_2)\bigr) \\
 & = &
 iZ(0)\Biggl(\frac{1}{\bigl|\Aut\Gamma_1\bigr|}W(\Gamma_1)
   +  \frac{1}{\bigl|\Aut\Gamma_2\bigr|}W(\Gamma_2)\Biggr)
\end{eqnarray*}
where $\Aut\Gamma_1$ and $\Aut\Gamma_2$ are the
automorphism groups of the graphs. We think of a trivalent graph as a
collection of half-edges partitioned into two element subsets (edges)
and three element subsets (vertices). Then in this context, an
automorphism of the graph means a permutation of the half-edges which
preserves both the partition into edges and the partition into
vertices (this differs from the usual definition). For example,
$\Aut\Gamma_1$ is generated by the involutions
$$\involutionI,\qquad\involutionII,\qquad\mbox{and}\qquad\involutionIII,$$
whereas $\Aut\Gamma_2$ is the semi-direct product of the group
generated by the involution
$$\involutionIV$$
and the symmetric group $S_3$ permuting the three edges.

Using the same arguments, we deduce that the $2m^{\mathrm{th}}$ term
looks like
$$\frac{(-1)^m}{(2m)!6^{2m}} V_{jkl}\ldots
V_{abc}\,i^{3m}Z(0)\bigl(Q^{jk}\ldots Q^{bc}+\cdots \bigr)$$
with $2m$ copies of $V$, and the last factor being the sum over all
ways of dividing a set of $3\times 2m$ indices into pairs, thus a sum
of $\frac{1}{(3m)!}\binom{6m}{\,2,\ldots,2\,}$ terms. Each trivalent
graph $\Gamma$ with $2m$ vertices gives us a weight $W(\Gamma)$ as
before, which we claim occurs precisely
$\smash{\frac{(2m)!6^{2m}}{|\Aut\Gamma|}}$ times in the sum above. To
understand why, it helps to separate the trivalent vertices from the
edges of the graph.
$$2m\left\{\begin{array}{cc}
   V & \vertex \\
     & \vdots \\
   V & \vertex \end{array}\right.\qquad\left.\begin{array}{cc}
   \edge & Q^{-1} \\
   \edge & Q^{-1} \\ 	
   \vdots & \\
   \edge & Q^{-1} \end{array}\right\}3m$$
Then $\Gamma$ is given by gluing together the legs in some specific
way. Note that if the legs are labeled with indices, then gluing
corresponds to contracting indices. Now observe that the semi-direct
product of groups $S_{2m}\ltimes (S_3)^{2m}$ acts on the trivalent
vertices. Under this action, we still get $\Gamma$ but the dummy
indices change, unless the element of the group induces an
automorphism of $\Gamma$. Thus the orbit has size
$\smash{\frac{(2m)!6^{2m}}{|\Aut\Gamma|}}$, corresponding to the
number of times $W(\Gamma)$ occurs in the sum above.

Thus the $2m^{\mathrm{th}}$ term is
$$i^mZ(0)\sum_{(m+1)-\mathrm{loop}\atop\mathrm{graphs}}\frac{1}{|\Aut\Gamma|}W(\Gamma)$$
(trivalent graphs with $2m$ vertices are usually referred to as
$(m+1)$--loop graphs), and
$$Z_k=\sum_{m=0}^{\infty}i^mZ(0)k^{-\frac{N}{2}-m}
  \sum_{(m+1)-\mathrm{loop}\atop\mathrm{graphs}}\frac{1}{|\Aut\Gamma|}W(\Gamma).$$
This sum over graphs includes disconnected graphs. However, the
multiplicativity of $W(\Gamma)$ and the simple behaviour of
$\Aut\Gamma$ under disjoint union readily leads to the
following formula:
$$\log\left(\frac{Z_k}{Z(0)}\right) =
  -\frac{N}{2}\log k
  +\sum_{m=0}^{\infty}i^mk^{-m}
  \sum_{\mathrm{connected}\atop{(m+1)-\mathrm{loop}\atop\mathrm{graphs}}}
  \frac{1}{|\Aut\Gamma|}W(\Gamma)$$

\section{Degenerate quadratic forms}

Before returning to the infinite-dimensional case, we need to consider
how to modify our formulae in the case that the quadratic form is
degenerate. Recall that the Gaussian integral for a non-degenerate
form $Q$ is given by
$$\int_{\mathbb{R}^N}d^Nxe^{ik\frac{1}{2}Q(x,x)}=\left(\frac{2\pi}{k}\right)^{\frac{N}{2}}\frac{e^{\pi
i(\mathrm{sign}Q)/4}}{|\det Q|^{\frac{1}{2}}}.$$
Now suppose that a compact group $K$ acts on $\mathbb{R}^N$ preserving
$Q$. We can write
$$\int_{\mathbb{R}^N}d^Nxe^{ik\frac{1}{2}Q(x,x)}=\int_{\mathbb{R}^N/K}d\mu
e^{ik\frac{1}{2}Q(x,x)}\Vol \mathcal{O}_x$$
where $d\mu$ is the induced measure on the space of orbits
$\mathrm{R}^N/K$, and $\Vol \mathcal{O}_x$ is the volume of the
orbit through $x\in\mathbb{R}^N$. More precisely, the integral over
$\mathrm{R}^N/K$ can be thought of as an integral over a slice for the
$K$--action, that is, a submanifold of $\mathrm{R}^N$ intersecting the
generic orbit in exactly one point. The action induces a map
$$\rho_x\co K\rightarrow\mathcal{O}_x\subset\mathbb{R}^N$$
taking $g\in K$ to $g(x)\in\mathcal{O}_x$. Denote the derivative of
this map by 
$$B:=d\rho_x\co\mathfrak{k}\rightarrow
T_x\mathcal{O}_x\subset T_x\mathbb{R}^N=\mathbb{R}^N$$
where $\mathfrak{k}$ is the Lie algebra of $K$. Then
\begin{eqnarray*}
\Vol \mathcal{O}_x & = & \int_{\mathcal{O}_x}d\eta^{\prime} \\
   & = & \int_Kd\eta|\det B| \\
   & = & \Vol K|\det B^*B|^{\frac{1}{2}}
\end{eqnarray*}
where $d\eta$ is the Haar measure on $K$, $d\eta^{\prime}$ is the
induced measure on $\mathcal{O}_x$, and $B^*$ is the adjoint of
$B$. Since $\rho_x$ is $K$--equivariant, $\det B$ is constant on $K$,
justifying the last equality. We usually refer to $|\det
B^*B|^{\frac{1}{2}}$ as the Jacobian determinant. Summarizing, we have
\begin{eqnarray}
\label{invariant}
\int_{\mathbb{R}^N}d^Nxe^{ik\frac{1}{2}Q(x,x)} & = & \Vol K\int_{\mathbb{R}^N/K}d\mu
e^{ik\frac{1}{2}Q(x,x)}|\det B^*B|^{\frac{1}{2}}.
\end{eqnarray}

\begin{exm}
The usual $S^1$--action on $\mathbb{R}^2$ preserves the standard
quadratic form ${x^2+y^2}$. In this case, $B_{(x,y)}$ maps $\theta$ to
$(-\theta y,\theta x)$ and the adjoint $B^*_{(x,y)}$ maps $(X,Y)$ to
${-Xy+Yx}$, so 
$$|\det B^*B|^{\frac{1}{2}}=(x^2+y^2)^{1/2}=r$$ 
and thus
$$\int_{\mathbb{R}^2}dxdye^{-(x^2+y^2)}=2\pi\int_{r=0}^{\infty}dr e^{-r^2}r=\pi.$$
\end{exm}

In the example above, $Q$ is still non-degenerate. However, if the
action of the group $K$ is affine, then this forces the form $Q$ to be
degenerate in the direction of the orbits. The author has been unable
to find such an example with $K$ compact, though in
infinite-dimensions the action of the gauge group on the space of
connections is affine. Now when $K$ is non-compact, a $K$--invariant
function will be non-integrable. However, the integral on the
right-hand side of \eqref{invariant} can still be finite, and we will
use this as a {\em definition} of the left-hand side (we will also
ignore the infinite $\Vol K$ factor in this case).

We assume that $Q$ induces a non-degenerate form on the space of
orbits. Thus we can reduce the original ill-defined integral involving
a degenerate quadratic form to one over the space of orbits, with the
inclusion of a factor proportional to the Jacobian determinant. This
additional factor can be dealt with via the Faddeev--Popov method,
which we will describe shortly, but first we will compute the
one-loop term directly in the infinite-dimensional case.

\section{The one-loop term}

In this section we will calculate the one-loop term in perturbative
Chern--Simons theory following Atiyah~\cite{atiyahNineZero}. Recall that the
contribution to the Chern--Simons path integral coming from connections
$A=\alpha+\beta$ in a neighbourhood of the flat connection $\alpha$ is
$$\int_{\beta\in\Omega^1(M)\otimes\mathfrak{g}}d\beta
e^{ikS(\alpha)}e^{\frac{ik}{4\pi}\int_M\Tr(\beta\wedge
d_{\alpha}\beta+\frac{2}{3}\beta\wedge\beta\wedge\beta)}.$$
In this section we will ignore the cubic term in $\beta$ and just
focus on the infinite-dimensional Gaussian integral coming from the
quadratic term, which is known as the one-loop term. In other words,
we will evaluate
$$\int_{\beta\in\Omega^1(M)\otimes\mathfrak{g}}d\beta
e^{\frac{ik}{4\pi}\int_M\Tr(\beta\wedge d_{\alpha}\beta)}.$$

First observe that after choosing a Riemannian metric on $M$, we can
define a pairing on $\Omega^1(M)\otimes\mathfrak{g}$ by
$$\langle\beta_1,\beta_2\rangle:=-\frac{1}{2\pi}\int_M\Tr(\beta_1\wedge
*\beta_2)$$
where $*$ is the Hodge star operator acting on the one-form part of
$\beta_2$ (which of course depends on the choice of metric). Recall
that for $X$ and $Y\in\mathfrak{su}(n)$, $-\Tr(XY)$ is a
positive-definite invariant bilinear form, which explains why we have
included a minus sign in the above pairing. The quadratic form $Q$ is
then given by the self-adjoint operator $-*d_{\alpha}$ on
$\Omega^1(M)\otimes\mathfrak{g}$, as
$${\textstyle\frac{1}{2}}ik\bigl\langle\beta,Q(\beta)\bigr\rangle=\frac{ik}{4\pi}\int_M\Tr(\beta\wedge
d_{\alpha}\beta).$$
This is completely analogous to the finite-dimensional case, where the
quadratic form is given by a matrix combined with the standard pairing
on $\R^N$.

Now the operator $Q=-*d_{\alpha}$ is degenerate on
$\Omega^1(M)\otimes\mathfrak{g}$, and the analogue of
$B\co\mathfrak{k}\rightarrow T_x\mathcal{O}_x\subset\mathbb{R}^N$ is
the map coming from infinitesimal gauge transformations
$$d_{\alpha}\co\Omega^0(M)\otimes\mathfrak{g}\rightarrow\Omega^1(M)\otimes\mathfrak{g}.$$
In other words, $Q$ is degenerate on the image
$d_{\alpha}\Omega^0$ (we use $\Omega^j$ as an abbreviation for
$\Omega^j(M)\otimes\mathfrak{g}$). If we assume that the flat
connection $\alpha$ is irreducible and isolated, then the twisted
de Rham complex
$$\Omega^0(M)\otimes\mathfrak{g}\stackrel{d_{\alpha}}{\longrightarrow}\Omega^1(M)\otimes\mathfrak{g}\stackrel{d_{\alpha}}{\longrightarrow}\Omega^2(M)\otimes\mathfrak{g}\stackrel{d_{\alpha}}{\longrightarrow}\Omega^3(M)\otimes\mathfrak{g}$$ 
is exact. Specifically, a flat connection gives a representation of
$\pi_1(M)$ which is irreducible if and only if
$\mathrm{H}^0(M,\mathfrak{g})$ vanishes, and
$\mathrm{H}^1(M,\mathfrak{g})$ is the space of first-order
deformations of the connection, which vanishes if the connection is
isolated. The remaining cohomology groups $\mathrm{H}^2$ and
$\mathrm{H}^3$ are dual to $\mathrm{H}^1$ and $\mathrm{H}^0$
respectively.

According to the previous section, our Gaussian integral is defined as
$$\int_{\beta\in\Omega^1/d_{\alpha}\Omega^0}d\beta
e^{ik\frac{1}{2}\langle\beta,Q^{\prime}(\beta)\rangle}|\det B^*B|^{\frac{1}{2}}.$$
The quadratic form
$Q^{\prime}$ is given by the reduction of $Q$ on the quotient
${\Omega^1/d_{\alpha}\Omega^0}$; note that it is non-degenerate by
exactness of the de Rham complex.

Define an operator $P$ on 
$$\Omega^0\oplus\Omega^1=\Omega^0\oplus d_{\alpha}\Omega^0\oplus(\Omega^1/d_{\alpha}\Omega^0)$$ 
by $(d_{\alpha},d_{\alpha}^*,-*d_{\alpha})$ where
$d^*_{\alpha}:=-*d_{\alpha}*$ is the adjoint of $d_{\alpha}$ on
one-forms. Writing this as a matrix (with respect to the three
factors)
$$P=\left(\begin{array}{ccc}
   0 & B^* & 0 \\
   B & 0 & 0 \\
   0 & 0 & Q^{\prime}
   \end{array}\right)$$
we see that
$$|\det P|=|\det B^*B||\det Q^{\prime}|.$$
Moreover,
\begin{eqnarray*}
B^*B & = & \Delta^0_{\alpha}\hspace*{20mm}(\mbox{the Laplacian on }\Omega^0(M)\otimes\mathfrak{g}) \\
P^2  & = & \Delta^0_{\alpha}\oplus\Delta^1_{\alpha}\hspace*{10mm}(\mbox{the
Laplacian on }\Omega^0(M)\otimes\mathfrak{g}\oplus\Omega^1(M)\otimes\mathfrak{g}).
\end{eqnarray*}
Note that although $Q$ is degenerate on $\Omega^1$, we can extend it
to a non-degenerate quadratic form $P$ on the larger space
$\Omega^0\oplus\Omega^1$. (This partly explains the addition of the
ghost field $c\in\Omega^0$ in the Faddeev--Popov method -- see the next
section.)

The Jacobian determinant $|\det
B^*B|^{\frac{1}{2}}=(\det\Delta^0_{\alpha})^{\frac{1}{2}}$ is
independent of $\beta$, while when we integrate the Gaussian we pick
up a factor of $|\det Q^{\prime}|^{-\frac{1}{2}}$. This combination
can be rewritten
\begin{eqnarray*}
\frac{|\det B^*B|^{\frac{1}{2}}}{|\det Q^{\prime}|^{\frac{1}{2}}} & =
& \frac{|\det B^*B|}{|\det P|^{\frac{1}{2}}} \\
   & = &
\frac{\bigl(\det\Delta^0_{\alpha}\bigr)^{\frac{3}{4}}}
     {\bigl(\det\Delta^1_{\alpha}\bigr)^{\frac{1}{4}}}.
\end{eqnarray*}
For the determinants of the Laplacians, we need to calculate the
product of their eigenvalues, which can be done using zeta-function
regularization. This gives the one-loop term, up to a phase factor and
an overall constant.

Bearing in mind that the Chern--Simons partition function is supposed
to represent a topological invariant of $M$, one might wonder about
the appearance of a Riemannian metric in our calculations above. While
it is true that the Laplacians (and their determinants) depend on the
choice of metric, Ray and Singer showed that the ratio
$$\frac{\bigl(\det\Delta^0_{\alpha}\bigr)^{\frac{3}{2}}}
       {\bigl(\det\Delta^1_{\alpha}\bigr)^{\frac{1}{2}}}$$
is independent of this choice~\cite{rsSevenOne}. They showed this by
explicitly calculating the variation under a change of metric, and
this quantity became known as the Ray--Singer analytic torsion. It was
later shown that it equals the Reidemeister torsion, which is defined
combinatorially. In particular, our factor above is the square root of
the Ray--Singer torsion, and is therefore a topological
invariant.

\section{The Faddeev--Popov method}

There is another way of dealing with the Jacobian determinant, which
is necessary to extend our results to the higher order terms in the
perturbative expansion. This is the method used by Axelrod and
Singer~\cite{asNineTwo} and Bar-Natan~\cite{bar-natanNineOne,bar-natanNineFive}, and  
involves introducing additional fields and integrating over them in
the path integral, as we now explain.

To begin with, we want to integrate over the space of orbits of the
gauge group acting on $\Omega^1(M)\otimes\mathfrak{g}$. One way to
achieve this is to find a function $F$ whose level sets meet each
orbit transversely in a single point, and then integrate over
$F^{-1}(0)$. Of course, we must also include the Jacobian determinant,
which represents the volume of the orbit. Equivalently, we could
include a $\delta$--function term $\delta(F(\beta))$ in the path
integral. Up to a factor, this can be written as
$$\delta(F(\beta))\propto \int d\phi e^{i\langle F(\beta),\phi\rangle}$$
by using the Fourier representation of the $\delta$--function. The term
$F(\beta)\phi$ can then be added to the Lagrangian, and $\phi$ becomes
an additional field to integrate over. In the Chern--Simons case, $F$
comes from the gauge-fixing condition $d^*_{\alpha}\beta=0$, and
$\phi$ should be a $\mathfrak{g}$--valued function.

The Jacobian determinant can also be written in terms of additional
fields. For this we must introduce ghost fields $c$ and $\overc$,
which lie in $\Omega^0(M)\otimes\mathfrak{g}$ but are Fermionic, that is,
anti-commutative (whereas $\beta$ and $\phi$ are Bosonic). Then, up to
a factor,
$$|\det B^*B|^{\frac{1}{2}}=
  \bigl(\det\Delta_{\alpha}^0\bigr)^{\frac{1}{2}}\propto \int
dc\,d\overc\, e^{ik\frac{1}{2}\langle\overc,\Delta_{\alpha}^0c\rangle},$$ 
which is the analogue of the Gaussian in the theory of
anti-commutative integration.

So our new path integral looks like
$$\int d\beta\, d\phi\, dc\,d\overc\, e^{ikS(\beta,\phi,c,\overc)}$$
where after the addition of gauge-fixing and ghost terms the action is
given by
$$S(\beta,\phi,c,\overc)=\frac{1}{4\pi}\int_M\Tr\bigl(\beta\wedge
d_{\alpha}\beta+{\textstyle\frac{2}{3}}\beta\wedge\beta\wedge\beta
+d^*_{\alpha}\beta\wedge*\phi+\overc\wedge *\Delta_{\alpha}^0c\bigr).$$
One might be concerned about the appearance of metric dependent terms
in the new Lagrangian, such as the $*$--operator and Laplacian. To show
that the path integral is independent of the choice of metric we
introduce the BRST operator. In this subsection only, $Q$ will denote
the BRST operator, which should not be confused with the quadratic
form $Q$ defined earlier. It is an odd (interchanges Bosonic and
Fermionic fields) derivation on the space of functionals in $\beta$,
$\phi$, $c$, and $\overc$ defined by
\begin{eqnarray*}
Q\beta & = & -(d_{\alpha}+\beta)c \\
Q\phi & = & 0 \\
Q\overc & = & \phi \\
Qc & = & {\textstyle\frac{1}{2}}[c,c].
\end{eqnarray*}
When we deform the metric, only the gauge-fixing and ghost terms
vary. However, one can show that this variation is
$Q$--exact. Moreover, $Q$ is like a vector field with zero divergence,
so the integral of a $Q$--exact term must vanish. See Axelrod and
Singer~\cite{asNineTwo} or Bar-Natan~\cite{bar-natanNineFive} for the details.

The introduction of the BRST operator can perhaps be understood as an
algebraic method of forming a quotient, though as far as the author
knows, this idea has not yet been mathematically formalized\footnote{Note
added in proof: for a homological algebraic interpretation due to
Jim Stasheff see \cite{stasheffTwo,stasheffThree,stasheffOne}.}. In our
case, we are interested in the quotient of the space of connections by
gauge transformations. Note that $Q$ extends the infinitesimal gauge
action on $\beta$ to the other fields $\phi$, $c$, and $\overc$. It
is straightforward to check that $Q^2=0$; essentially what we are
doing above is taking $Q$--cohomology. Only the metric invariant part
$$\frac{1}{4\pi}\int_M\Tr\bigl(\beta\wedge
d_{\alpha}\beta+{\textstyle\frac{2}{3}}\beta\wedge\beta\wedge\beta\bigr)$$
survives to give an element of $\mathrm{H}^0(Q)$, which is our
algebraic representation of the space of functionals on the
quotient. The BRST approach to gauge theory is described in Witten's
third lecture in~\cite{wittenNineNineii}; Chapter 7 of Bar-Natan's
thesis~\cite{bar-natanNineOne} is another useful reference.

It should also be noted that a more natural way to view the above path
integral is in its superspace formulation. However, since we wish to
write the terms of the perturbative expansion as configuration space
integrals we must stick to the language of differential forms. Both of
these approaches are described by Axelrod and Singer~\cite{asNineTwo}.

\section{Higher order terms}

We will now give a rough outline of how one constructs the higher
order terms in the perturbative expansion of the Chern--Simons
partition function. These terms will be indexed by trivalent graphs,
and formally the series will look the same as the finite-dimensional
case, once we have properly interpreted the quadratic and cubic forms
on $\Omega^1(M)\otimes\mathfrak{g}$. Actually, the series has only
been fully computed in the case that $\alpha$ is the trivial
connection, in which case the twisted de Rham complex decouples 
$$\Omega^0(M)\otimes\mathfrak{g}\stackrel{d\otimes\mathrm{Id}}{\longrightarrow}\Omega^1(M)\otimes\mathfrak{g}\stackrel{d\otimes\mathrm{Id}}{\longrightarrow}\Omega^2(M)\otimes\mathfrak{g}\stackrel{d\otimes\mathrm{Id}}{\longrightarrow}\Omega^3(M)\otimes\mathfrak{g}$$  
where $d$ is the usual exterior derivative on forms on $M$, and
$\mathrm{Id}$ is the identity on $\mathfrak{g}$. We will assume
$\alpha$ is trivial for the remainder of this article; of course if
$M$ is Euclidean space then the only flat connection is the trivial
connection.

We saw that the quadratic form $Q$ is given by the 
self-adjoint operator $-*d_{\alpha}$, when combined with the pairing
$$\langle\beta_1,\beta_2\rangle:=-\frac{1}{2\pi}\int_M\Tr(\beta_1\wedge
*\beta_2)$$
on $\Omega^1(M)\otimes\mathfrak{g}$. When $\alpha$ is the trivial
connection this decouples into $-*d$ (combined with the induced
inner product) on $\Omega^1(M)$ and the invariant quadratic form $B$
on $\mathfrak{g}$. Recall that to calculate weights on trivalent graphs
we needed the inverse $Q^{-1}$. As far as $\mathfrak{g}$ is concerned,
this is just the inverse $B^{-1}$ of the invariant quadratic form. 

On $\Omega^1(M)$ we need to calculate the inverse of the elliptic
operator $-*d$. Strictly speaking, this inverse does not exist as
$-*d$ has non-trivial kernel. However, the kernel is the image of the
infinitesimal gauge action, and is taken care of by the Faddeev--Popov
method described in the last section. After adding gauge-fixing terms,
we find that the new Lagrangian has a quadratic term which looks like
$$ik{\textstyle\frac{1}{2}}\langle \beta+\phi, L(\beta+\phi)\rangle$$
where $L:=-(*d+d*)$. Moreover, $L^2$ is the Laplacian (on functions
and one-forms), whose inverse is given in terms of a Green's
function. For example, in Euclidean space $\mathbb{R}^3$ the solution
to $\Delta u=w$ is given by
$$u(x)=\frac{1}{4\pi}\int_{\mathbb{R}^3}dy\frac{w(y)}{|x-y|}.$$
The Green's function is  
$$G(x,y)=\frac{1}{4\pi |x-y|}$$
and the inverse of $L$ itself is given by 
$$L\circ G=\frac{\epsilon^{ijk}(x_i-y_i)dx_j\wedge
  dx_k}{4\pi |x-y|^3}$$
(using Einstein summation), meaning that this is the kernel of an
integral operator representing $L^{-1}$.

The cubic term $V$ is much simpler. Recall that if we write
$\beta=\zeta^IE_I$ where $\zeta^I$ are one-forms and $\{E_I\}$ 
is a basis for $\mathfrak{g}$, then
$$\Tr(\beta\wedge\beta\wedge\beta)=
  {\textstyle\frac{1}{2}}c_{IJK}\zeta^I\wedge\zeta^J\wedge\zeta^K$$ 
where $c_{IJK}$ comes from the structure constants of
$\mathfrak{g}$. So the cubic term decouples into simply wedging the
three one-forms together (and then integrating over $M$) on
$\Omega^1(M)$ and $c_{IJK}$ on $\mathfrak{g}$.

The Faddeev--Popov method also introduces ghost fields, so there is an
additional propagator and interaction coming from ghosts. These are
usually drawn with dashed lines in our Feynman diagrams, whereas the
propagator and interaction described above are drawn with solid
lines. See Bar-Natan~\cite{bar-natanNineOne,bar-natanNineFive} for more details.

Let $\Gamma$ be a trivalent graph. Since both the quadratic and cubic
terms decouple, the weight $W(\Gamma)$ will become a product
\begin{eqnarray}
\label{product}
W(\Gamma) & = & I_{\Gamma}(M)\times b_{\mathfrak{g}}(\Gamma)
\end{eqnarray}
where $I_{\Gamma}(M)$ does not depend on the Lie algebra
$\mathfrak{g}$ and $b_{\mathfrak{g}}(\Gamma)$ does not depend on the
three-manifold $M$. The fact that $\alpha$ is the trivial connection
is crucial for this decoupling, as otherwise the quadratic form would
not decouple.

\begin{rem}
For each finite-dimensional Lie algebra $\mathfrak{g}$,
$b_{\mathfrak{g}}$ is a weight system on trivalent graphs. It is the
same as the weight system arising in the theory of Vassiliev
invariants~\cite{bar-natanNineFiveii}; for example
$$b_{\mathfrak{g}}(\Gamma_2):=c_{IJK}c_{ABC}B^{IA}B^{JB}B^{KC}.$$
\end{rem}

Notice that the cubic terms on $\Omega^1(M)$ and $\mathfrak{g}$ are
both skew-symmetric. Consequently we really need to introduce the
notion of an orientation on $\Gamma$ (an equivalence class of cyclic
orderings of the legs at each trivalent vertex), so that
$I_{\Gamma}(M)$ and $b_{\mathfrak{g}}(\Gamma)$ make sense. The
skew-symmetry cancels out in $W(\Gamma)$: if we reverse the
orientation on $\Gamma$, both $I_{\Gamma}(M)$ and
$b_{\mathfrak{g}}(\Gamma)$ will change sign but $W(\Gamma)$ will
remain the same. 

One consequence of this skew-symmetry of vertices is that many weights
automatically vanish. For example, write $\Gamma_1$ and $\widetilde{\Gamma}_1$
for the graph
$$\dumbell$$
with its two orientations. Then $\Gamma_1$ and $\widetilde{\Gamma}_1$ are
actually isomorphic as oriented graphs, so
$$b_{\mathfrak{g}}(\Gamma_1)=b_{\mathfrak{g}}(\widetilde{\Gamma}_1)=-b_{\mathfrak{g}}(\Gamma_1)$$
must vanish.

\section{Configuration space integrals}

We will describe in more detail the terms $I_{\Gamma}(M)$ arising
from \eqref{product}. In our finite-dimensional model we summed over
indices when evaluating $W(\Gamma)$, as
in \eqref{summation_convention}, and this still works fine for
the Lie algebra factor $b_{\mathfrak{g}}(\Gamma)$ since $\mathfrak{g}$
is finite-dimensional. However, on the infinite-dimensional part
summation is replaced by integration: heuristically, one can regard
the $\delta$--functions $\{\delta(x-y)\}_{y\in M}$ as a basis for
functions on $M$, so that summation over elements in the basis becomes 
$\int_{y\in M}$ (and something similar happens for one-forms). The
term $I_{\Gamma_2}(M)$, for example, therefore looks like
$$\int_Mdx\int_Mdy K(x,y)\wedge K(x,y)\wedge K(x,y).$$
In this formula, $K(x,y)$ is the Green's function for $L^{-1}$ on
$M$. We can only write this out explicitly when $M$ is Euclidean
space, since we have an explicit formula for the Green's function
on $\mathbb{R}^3$.

For a graph $\Gamma$ with $m$ vertices, we get an integral over
$M^m$. Now the Green's function $K(x,y)$ will have a singularity at
$x=y$, so the integrand only makes sense on the configuration space
$$C^0_m(M):=\{(x^{(1)},\ldots,x^{(m)})\in M^m~|~x^{(i)}\neq x^{(j)}\mbox{ for }i\neq
j\}.$$
However, by analyzing the severity of the singularity Axelrod and
Singer~\cite{asNineTwo} showed that the integrals on $M^m$ are actually
finite. One approach here is to extend the integrand to the
compactification of $C^0_m(M)$ due to Fulton and
MacPherson~\cite{fmNineFour}, which is a manifold with boundary and corners,
and then show that this extension is finite.

It is important to note that for a fixed trivalent graph $\Gamma$,
$I_{\Gamma}$ is {\em not\/} a topological invariant of
three-manifolds: it varies as we deform the metric on $M$. The
variation can be understood as an integral over the boundary of the
Fulton--MacPherson compactification of $C^0_m(M)$. However, when we
combine the terms $I_{\Gamma}$ in the sum
\begin{eqnarray}
\label{combination}
\sum_{(m+1)-\mathrm{loop}\atop\mathrm{graphs}}\frac{1}{|\Aut\Gamma|}I_{\Gamma}(M)\times
b_{\mathfrak{g}}(\Gamma)
\end{eqnarray}
the variations should cancel out. In fact there are two kinds of
variation, and variations of the first kind cancel due to a special
property of the weights $b_{\mathfrak{g}}(\Gamma)$, namely that they
satisfy the IHX relations. Variations of the second kind are known as
``anomalous'', and are more difficult to deal with. If we choose a
framing of $M$, then we can correct the anomaly with a counter-term
and thereby construct a genuine invariant of $M$ (see Sections 5 and 6
of Axelrod and Singer~\cite{asNineTwo}).

\section{Invariants of knots and links}

Finally let us say a few words about the generalization to knot and
link invariants. Given an oriented knot $\mathcal{K}$ in $M$, we can
include an additional term in the Chern--Simons path integral given by
the monodromy of the connection $A$ around $\mathcal{K}$. The
monodromy would take values in $\mathfrak{g}$, so we evaluate its
trace in a representation $V$ of $G$. This gives the Wilson
loop
$$\mathcal{W}_V(\mathcal{K}):=\Tr_V\mathrm{Pexp}\int_{\mathcal{K}\subset
  M}A.$$
In this expression, the path ordered exponential $\mathrm{Pexp}$ can
be represented by a sum of iterated integrals 
\begin{eqnarray}
\label{iterated}
\sum_{k=0}^{\infty}\mathop{\idotsint}_{0\leq t_1\leq t_2\leq\ldots\leq
  t_k\leq 1}A(\mathcal{K}(t_1))\otimes\cdots\otimes
A(\mathcal{K}(t_k))
\end{eqnarray}
where we have used $\mathcal{K}\co S^1\rightarrow M$ to denote a
parametrization of the knot. We then represent the
$\mathfrak{g}$--valued $A(\mathcal{K}(t_i))$ in $\mathrm{End}(V)$,
multiply them, and take the trace. The expectation value of the Wilson
loop
$$Z_k(M;\mathcal{K})=\int_{\mathcal{A}}\mathcal{D}Ae^{ikS(A)}\mathcal{W}_V(\mathcal{K})$$
should be an invariant of the knot. For links one includes several
Wilson loops, one for each component of the link.

The expectation value $Z_k(M;\mathcal{K})$ admits a perturbative
expansion just like $Z_k(M)$, and once again each term splits into the 
product of a weight system (coming from the Lie algebra $\mathfrak{g}$
and representation $V$) and a configuration space integral. The
Feynman graphs consist of unitrivalent graphs, whose univalent
vertices lie on a circle $S^1$ (or collection of $l$ circles in the
case of an $l$--component link). The interactions at the univalent
vertices come from the terms $A(\mathcal{K}(t_i))$
in \eqref{iterated}. The configuration space looks like
$$C^0_{a,b}(M;\mathcal{K}):=\bigl\{(x^{(1)}\!\!,\ldots,x^{(a)}\!\!,x^{(a+1)}\!\!,
  \ldots,x^{(a+b)})\in
  M^a\times\mathcal{K}^b\,|\,x^{(i)}\neq x^{(j)}\mbox{\,for\,}i\neq j\bigr\}$$
and the integral can be replaced by an integral over a suitable
compactification of $C^0_{a,b}(M;\mathcal{K})$ (see
Thurston~\cite{thurstonNineFive}).

A great deal of work has gone into understanding these configuration
space integrals when $M$ is Euclidean space, since we then have an
explicit formula for the Green's function of $L^{-1}$ (see
Bar-Natan~\cite{bar-natanNineOne,bar-natanNineFive}, Bott and Taubes~\cite{btNineFour},
Altschuler and Freidel~\cite{afNineSeven}, and
Thurston~\cite{thurstonNineFive}). The simplest example comes from the graph
given by a single chord connecting two circles, in the case of a two
component link $\mathcal{L}$ (it looks like $\Gamma_1$, but with the
circles corresponding to the two components of the link). Then 
$$C^0_{0,\{1,1\}}(M;\mathcal{L})\cong\mathcal{L}_1\times\mathcal{L}_2\cong
S^1\times S^1$$ 
is already compact (we write $\{1,1\}$ to indicate that there is a
point on each component of the link). The integral is given by
$$\int_{\mathcal{L}_1}dx_j\int_{\mathcal{L}_2}dy_k\frac{\epsilon^{ijk}(x_i-y_i)}{4\pi|x-y|^3}$$
and it computes the degree of the map from
$\mathcal{L}_1\times\mathcal{L}_2\cong S^1\times S^1$ to $S^2$ which
takes $(x,y)$ to $\frac{x-y}{|x-y|}$. It is known as the Gauss linking
number of $\mathcal{L}_1$ and $\mathcal{L}_2$.

The Gauss linking number is a genuine invariant, but for a general
configuration space integral there are two kinds of variation when we
isotope the knot or link. It is easy to show that, when combined in a
sum like \eqref{combination}, the variations of the first kind cancel
out. However, the second kind of variation is more difficult to deal
with, and once again leads to an ``anomaly''. 

For example, the first-order term
$I_{\Theta}(\mathbb{R}^3;\mathcal{K})$ for a knot $\mathcal{K}$ is
known as the writhe. It corresponds to the graph given by a circle
with a single chord (it looks like $\Gamma_2$, but with the circle
corresponding to the knot), and can be expressed as a Gauss integral
similar to the one above. This integral varies as we isotope the knot,
and can take all real values, so there is an anomaly. However, if the
knot is framed, that is, a trivialization of its normal bundle is given,
then we get a genuine invariant by adding the total torsion $\tau$
(see Pohl~\cite{pohlSixEight}). In fact, up to a factor,
$$I_{\Theta}(\mathbb{R}^3;\mathcal{K})+\tau$$
is the linking number of the knot with its parallel (the knot obtained
by using the framing to displace $\mathcal{K}$ slightly). 

It was shown independently by Altschuler and Freidel~\cite{afNineSeven} and
by Thurston~\cite{thurstonNineFive} that, up to a correction by the anomaly,
the perturbative series gives an invariant of knots and links in
$\mathbb{R}^3$. This extended earlier work by 
Bar-Natan~\cite{bar-natanNineOne,bar-natanNineFive} and Bott and
Taubes~\cite{btNineFour}, who had investigated the series up to the degree
two term. Moreover, the perturbative series was shown to be a universal
Vassiliev invariant (c.f.\ Bar-Natan~\cite{bar-natanNineFiveii}). It is
expected that it should agree with the universal Vassiliev invariant
of Kontsevich~\cite{kontsevichNineThree}, which would follow from the
vanishing of the anomaly in degree greater than one. This is still an
open conjecture, though recently Poirier~\cite{poirierZeroTwo} showed that
the anomaly does indeed vanish in degrees two to six (see also
Lescop~\cite{lescopZeroTwo}).

\bibliographystyle{gtart}
\bibliography{link}

\end{document}